\documentclass{article}
\usepackage{amsmath,amssymb,amsthm}
\usepackage{graphics}

\newtheorem{theorem}{Theorem}[section]

\newcommand{\dist}{\mathop{\rm dist}\nolimits}
\newcommand{\tr}{\mathop{\rm Tr}\nolimits}
\newcommand{\rk}{\mathop{\rm rank}\nolimits}

\title {Bounds for Codes by Semidefinite Programming }

\author {Oleg R. Musin \thanks{University of Texas at Brownsville. Email: omusin@gmail.com}}

\begin{document}
\date{}
\maketitle

\begin{abstract}
Delsarte's method and its extensions allow to consider the upper bound problem for codes in 2-point-homogeneous spaces as a linear programming problem with perhaps infinitely many variables, which are the distance distribution. We show that using as variables  power sums of distances this problem can be considered as a finite semidefinite programming  problem. This method allows to improve some linear programming upper bounds. In particular we obtain new bounds of one-sided kissing numbers. %We discuss a possibility to apply this approach for construction of  codes.
\end {abstract}

\section{Introduction}

The Delsarte linear programming  method is widely used for finding bounds for error-correcting codes, constant weight codes, spherical codes, sphere packings and other packing problems in 2-point-homogeneous spaces \cite{CS,Del1,Del2, Kab, Levm}.
Linear Programming (LP) is a special case of  Semidefinite Programming (SDP) which
deals with optimization problems over symmetric positive semidefinite matrix variables with linear cost function and linear constraints.

Recently, Schrijver \cite{Schr} using  SDP improved some upper bounds on binary codes. Even more recently, Schrijver's method has been adapted
 for  non-binary codes (Gijswijt,  Schrijver, and Tanaka \cite{GST}) and  for spherical codes  (Bachoc and Vallentin \cite{BV}). (See also our paper \cite{Mus4}.)
 In fact, this method using the stabilizer subgroup of the isometry group derives new positive semidefinite constraints which are stronger than the linear inequalities in the Delsarte linear programming method.

In this paper we define semidefinite programs whose optimal solutions give upper bounds for codes in
a 2-point-homogenous space and strengthen the LP bounds. Our approach is based on two ideas: to consider power sums of distances instead of the distance distribution, and to use the positive-semidefinite property of zonal spherical functions. In fact, this method is different from Schrijver's method and can be applied for any 2-point-homogenous space.

The paper is organized as follows: Section 2 reviews on the LP method for codes in 2-point-homogenous spaces. Section 3 recalls a theorem due to Sylvester. Section 4 introduces the $SDP_0$ bounds for codes and shows that these bounds are equivalent to the LP bounds. Section 5 considers examples of applications of $SDP_0$.
%for construction of  codes.
Section 6 shows that some recent extensions of Delsarte's method can be reformulated as SDP problems ($SDPA$). Section 7 extends the $SDPA$ bounds to subsets of a 2-point-homogeneous space and shows that some upper bounds for codes can be improved.

Our computations were performed by the program SDPA-M, version 2.0 (see http://sdpa.is.titech.ac.jp/SDPA-M/). After solving the SDP with SDPA-M we checked independently whether the solution satisfies the desired constraints. This can be done using rational arithmetic only.

\section{The linear programming bounds for codes}

The main goal of the present section is to give a brief introduction to zonal spherical functions and the linear programming bounds for codes. We are based on the simplified account, given by Sloane \cite[Chapter 9]{CS}, of the general machinery developed by Kabatiansky and Levenshtein \cite{Kab}.
Here we are using the same notations as in \cite{CS} and we hope that the reader will have no difficulty using Sloane's exposition to fill up details and proofs.

We say that a $G$-space ${\bf M}$ is a {\em 2-point-homogeneous space} if (i) ${\bf M}$ is a metric space with a distance  $\rho$ defined on it;
(ii) ${\bf M}$ is a set on which a group $G$ acts;  (iii) $\rho$ is strongly invariant under $G$, i.e. for $x,x',y,y' \in {\bf M}$ with $\rho(x,y)=\rho(x',y')$  there is an element $g\in G$ such that $g(x)=x'$ and $g(y)=y'$.

These assumptions are quite restrictive. In fact, if $G$ is infinite and ${\bf M}$ is a compact space, then  Wang \cite{Wang} has proved that ${\bf M}$ is a sphere; real, complex or quaternionic projective space; or the Cayley projective plane. However, the finite  2-point-homogeneous spaces have note yet been completely classified (see for the most important examples and references  \cite{CS}).

With any  2-point-homogeneous space ${\bf M}$ and an integer number $k\ge 0$ are associated the {\em zonal spherical function} $\Phi_k(t)$ such that $\{\Phi_k(t)\}_{k=0,1,2,\ldots}$  are orthogonal on
$T:=\{\tau(x,y): x, y\in {\bf M}\},$ where $\tau$ is the certain function in $\rho$ ( i.e
$\tau(x,y)=F(\rho(x,y))$) defined by ${\bf M}$.
For all continuous compact ${\bf M}$  and for all currently known finite cases: {\em $\Phi_k(t)$
 is a polynomial of degree } $k$. The normalization is given by the rule:
$\Phi_k(\tau_0)=1$, where $\tau_0:=\tau(x,x)$. Then  $\Phi_0(t)=1.$

Let us consider two examples.

\medskip

\noindent{\em Example 1:} ${\bf M}$ = Hamming space ${\bf F}_2^n$ with $\tau(x,y)=\rho(x,y)$=Hamming distance. Then $\tau_0=0$.  Here  $\Phi_k(t)$ is the Krawtchouk polynomial $K_k(t,n)$:
$$
\Phi_k(t)=K_k(t,n)={n\choose k}^{-1}\sum\limits_{j=0}^k {(-1)^j {t\choose j}{n-t\choose k-j}}.
$$

\noindent{\em Example 2} (the most extensively studied case): ${\bf M}$ = unit sphere ${\bf S}^{n-1}\subset {\bf R}^n$ with $\tau(x,y)=\cos{\rho(x,y)}=\langle x,y\rangle$, where  $\rho(x,y)$ is the angular distance between $x$ and $y$. In this case $\tau_0=1.$
The corresponding zonal spherical function $\Phi_k(t)$ is the Gegenbauer (or ultraspherical) polynomial $G_k^{(n)}(t)$.

There are many ways to define the Gegenbauer polynomial $G_k^{(n)}$ (see \cite{CS,Del2,Kab,Mus2,PZ,Scho}), which is a special case of  Jacobi polynomials
$P_k^{(\alpha,\beta)}$ with $\alpha=\beta=(n-3)/2$ and normalization  $G_k^{(n)}(1)=1$. $G_k^{(n)}$ are orthogonal on $[-1,1]$ with weight $(1-t^2)^{(n-3)/2}$. Also $G_k^{(n)}$ can be defined by the recurrence formula:
$$ G_0^{(n)}=1,\; G_1^{(n)}=t,\; \ldots,\; G_k^{(n)}=\frac {(2k+n-4)\,t\,G_{k-1}^{(n)}-(k-1)\,G_{k-2}^{(n)}} {k+n-3}.$$

\medskip

The main property for zonal spherical functions is called ``positive-definite degenerate kernels" or ``p.d.k"  \cite{CS}.
%(We think that for p.d.k. be better to use the term:``positive-semidefinite".)
This property first was discovered by Bochner \cite{Boc} (general spaces)  and independently for spherical spaces by Schoenberg \cite{Scho}.

Now we explain what the p.d.k. property means for finite subsets in ${\bf M}$.
%Let  $F_{i,j}^{(k)}=F_{i,j}^{(k)}(C):=\Phi_k(\tau(x_i,x_j))$.
\begin{theorem}[\cite{Boc,Scho,Kab}] Let ${\bf M}$ be a 2-point-homogeneous space.
Then for any integer $k\ge 0$ and for any finite  $C=\{x_i\}\subset {\bf M}$ the matrix $\left(\Phi_k(\tau(x_i,x_j))\right)$ is positive semidefinite.
\end{theorem}

This theorem implies the fact that plays a crucial role for the linear programming  bounds.
For any positive semidefinite matrix the sum of its entries is nonnegative. Then
\begin{theorem}[\cite{Del1,Del2,Kab,OdS}] For any finite  $C=\{x_i\}\subset {\bf M}$ we have
$$
\sum\limits_{i=1}^{|C|}\sum\limits_{j=1}^{|C|} {\Phi_k(\tau(x_i,x_j))}\ge 0.
$$
\end{theorem}

Since %$\tau(x_i,x_i)=\tau_0$, we see that
$\Phi_k(\tau(x_i,x_i))=\Phi_k(\tau_0)=1$, this theorem implies
$${\frac{1}{|C|} \sum\limits_{i,j:i\ne j} {\Phi_k(\tau(x_i,x_j))}\ge -1}, \; \, k=0, 1, 2, \ldots \, . \eqno (1)$$

Let $S$ be a fixed subset of $T$. We say that a finite subset $C\subset {\bf M}$ is an {\em $S$-code} if $\tau(x,y)\in S$ for all $x, y \in C, x\ne y.$ The largest cardinality $|C|$ of an $S$-code will be denoted by $A({\bf M},S)$.

%Coding theory and sphere packings are not interested in arbitrary $S$.
Usually, $S=T\cap[a,b]$  for some interval $[a,b]$.
For instance, if we take $S=T\cap[d,n]$ for Example 1, then an $S$-code is an  error-correcting code of minimal distance $d$.
%, and $A(M,S)=A(n,d).$
For Example 2, if $S=[-1,\cos{\theta}]$, an $S$-code (or a $\theta$-code) is a spherical code with minimal angular distance $\theta$.
% and $A(M,S)=A(n,\theta).$

 The {\em distance distribution} $\{\alpha_t\}$ of $C$ is defined by
$$
\alpha_t:=\frac{1}{|C|}(\mbox{number of ordered pairs }\, x, y\in C \, \mbox{ with }\, \tau(x,y)=t).
$$
We obviously have
$$
\alpha_{\tau_0}=1, \quad \alpha_t\ge 0, \; t\in T, \quad \sum\limits_{t\in T}{\alpha_t}=|C|. \eqno (2)
$$

$(1)$ and $(2)$ make it possible to regard the problem of bounding $A({\bf M},S)$ as a linear programming  problem:

\medskip

\noindent{\em Primal problem (LPP):} Choose a natural number $s$, a subset $\{\tau_1, \ldots, \tau_s\}$ of
  $S$, and real numbers $\alpha_{\tau_1}, \ldots, \alpha_{\tau_s}$ so as to
$$
\mbox{maximize }\; \alpha_{\tau_1} + \ldots + \alpha_{\tau_s}
%\sum\limits_{i=1}^s {\alpha_{\tau_i}}
$$
subject to
$$
\alpha_{\tau_i}\ge 0, \; i=1, \ldots, s, \quad \sum\limits_{i=1}^s {\alpha_{\tau_i}\Phi_k(\tau_i)}\ge -1, \; k=0,1, \ldots .
$$

This is a linear programming problem with perhaps infinitely many unknowns $\alpha_{t}$ and constraints $(1), (2)$.
If $C$ is an $S$-code then its distance distribution certainly satisfied the constraints $(1), (2)$. So if the maximal value of the sum $\alpha_{\tau_1} + \ldots + \alpha_{\tau_s}$ that can be attained is $A^*$, then $A({\bf M},S)\le 1+A^*.$ (The extra 1 arises because the term $\alpha_{\tau_0}=1$ doesn't occur in this sum.)

\medskip

\noindent{\em Dual problem (LPD):} Choose a natural number $N$ and real numbers  $f_1, \ldots, f_N$
 so as to
$$
\mbox{minimize }\; f_{1} + \ldots + f_N
%\sum\limits_{i=1}^s {\alpha_{\tau_i}}
$$
subject to
$$
f_k\ge 0, \; k=1, \ldots, N, \quad \sum\limits_{k=1}^N {f_k \, \Phi_k(t)}\le -1, \; \,  t\in S .
$$

Thus, we have
\begin{theorem} If $A^*$ is the optimal solution to either of the primal or dual problems, then
$A({\bf M},S)\le 1+A^*.$
\end{theorem}

Note that for the dual problem any feasible solution is an upper bound to the optimal solution of the primal problem. Using this (or using directly Theorem 2.2) it is not hard to prove:
\begin{theorem}[\cite{Del1,Del2,Kab,Lev2,OdS}] Suppose that $f$ is an nonnegative  linear combination of $\; \{\Phi_k\}_{k=0,1,\ldots,N} \; $ for some $N$, i.e.
$$
f(t)=f_0 \Phi_0(t)+f_1 \Phi_1(t)+\ldots +f_N \Phi_N(t), \; \mbox{ where } \; f_k\ge 0.
$$
If $f(t)\le 0\; $ for all $\; t\in S$ and $f_0>0$, then
$$
A({\bf M},S)\le \frac{f(\tau_0)}{f_0}.
$$
\end{theorem}

This theorem and its extensions have numerous  applications in coding theory and sphere packings (see \cite {CS}  for references, and for some recent achievements \cite{Levm,PZ,Mus,Mus13,Mus2,Mus3}). Since one of the most important steps in this approach was taken by Delsarte \cite {Del1}  it is
often  called Delsarte's linear programming bound or Delsarte's method.

\section{Sylvester's Theorem}

%Sylvester, Hermite and Sturm studied

One of the classical problems for polynomials is the problem of finding the number of positive roots of a real polynomial. There are  several  fundamental facts in this field (see \cite[Chapter 1]{Pras}). One of the powerful methods is called Hermite's method.  This method is
based on (Hankel's) quadratic forms which coefficients depend on power sums $s_k:=t_1^k+\ldots+t_n^k$, where $t_i$ are roots of a polynomial. Actually, this method gives an answer for the following question: Suppose we know for complex numbers $t_1, \ldots, t_n$ only its power sums $ s_k$, where  $ s_k$ are real numbers. How to determine the number of real and positive (or nonnegative) $t_i$?

%Let $f(t)$ be a real polynomial of degree $n$ with roots $t_1,\ldots,t_n$. %For simplicity, assume that $t_i\ne t_j$.
%For complex numbers $t_1, \ldots, t_n$  denote by $r(t_1,\ldots,t_n)$ the number of distinct real $t_i$, i.e. if $\{t_i\}=\{x_1,\ldots,x_m\}\cup \{y_i\}$, where $x_i$ are real numbers, and $y_i$ are not real, then $r=r(t_1,\ldots,t_n)$ is the minimal number such that $x_1, \ldots, x_m$ are real number
%$z_1, \ldots, z_r$ with multiplicities $m_1,\ldots,m_r, \; m_1+\ldots+m_r=m$.
%$s_k:=t_1^k+\ldots+t_n^k.$

\begin{theorem}[Sylvester, {\cite[Section 1.4]{Pras}}] $\\$
(i) Suppose that for complex numbers $t_1, \ldots, t_n$ all power sums $s_k$ are real numbers.
% Let $r(t_1,\ldots,t_n)$ be
Then the number of distinct real $t_i$
%Then  $r(t_1,\ldots,t_n)$
is equal to the signature of the quadratic form defined by the matrix
$$
R_n=
\left(
\begin{array}{cccc}
 s_0 &  s_1
&\ldots &  s_{n-1}\\
 s_1 &  s_2
&\ldots &  s_n\\
\vdots & \vdots
&\ddots & \vdots\\
 s_{n-1} &  s_n
&\ldots &  s_{2n-2}
\end{array}
\right).
$$

\noindent (ii) If  all $t_i$ are real, then  all of them are nonnegative if and only if the  matrix $F_n$ is positive semidefinite, where

$$
F_n=
\left(
\begin{array}{cccc}
s_1 &  s_2
&\ldots &  s_n\\
 s_2 &  s_3
&\ldots &  s_{n+1}\\
\vdots & \vdots
&\ddots & \vdots\\
 s_{n} &  s_{n+1}
&\ldots &  s_{2n-1}
\end{array}
\right).
$$

\end{theorem}

So we have that all $t_i$ are real and nonnegative if and only if $R_n$ and $F_n$ are positive semidefinite matrices ($R_n\succeq 0,\; F_n\succeq 0$), and the number of distinct $t_i$ equals $\rk(R_n)$.  It is clear that Sylvester's theorem  gives an answer for the same question when all $t_i$ are real and belong to some interval $[a,b]$.

%Let $s_k$ be normalized $\tilde s_k$, i.e. for given numbers $t_1, \ldots, t_n$

Now we introduce the following matrices:
$$
F_m^+(a):=
\left(
\begin{array}{cccc}
s_1-as_0 &  s_2-as_1
&\ldots &  s_m-as_{m-1}\\
s_2-as_1 & s_3-as_2
&\ldots &  s_{m+1}-as_m\\
\vdots & \vdots
&\ddots & \vdots\\
s_{m}-as_{m-1} &  s_{m+1}-as_m
&\ldots &  s_{2m-1}-as_{2m-2}
\end{array}
\right),
$$
$$
F_m^-(b):=
\left(
\begin{array}{cccc}
bs_0-s_1 &  bs_1-s_2
&\ldots &  bs_{m-1}-s_m\\
bs_1-s_2 & bs_2-s_3
&\ldots &  bs_{m}-s_{m+1}\\
\vdots & \vdots
&\ddots & \vdots\\
bs_{m-1}-s_{m}&  bs_{m}-s_{m+1}
&\ldots &  bs_{2m-2}-s_{2m-1}
\end{array}
\right),
$$
and
$$
H_m(a,b)=H_m(s_0,s_1\ldots,s_{2m-1},[a,b]):=
\left(
\begin{array}{cccc}
R_m & 0 & 0\\
0& F_m^+(a) &  0\\

0 & 0 & F_m^-(b)
\end{array}
\right).
$$

\begin{theorem} Consider real numbers $t_1,\ldots,t_n$ in $[a,b]$. Then for any natural number $m$ the matrix $H_m(a,b)$ is positive semidefinite. Moreover, for  $m\ge \rk(R_n)$ the converse holds, i.e. if $H_m(a,b)\succeq 0$, then  complex numbers  $t_1,\ldots,t_n$ with real $s_k$ are real numbers in $[a,b]$.% if and only if for   $m\ge r$:
%and $\rk(R_m)=r$.
\end{theorem}
\begin{proof} The proof is based on the classical Hermite's proof of Sylvester's theorem.

Let
$$
G({\bf u},{\bf v},{\bf w})=\sum\limits_{k=1}^n{(u_1+t_ku_2+\ldots+t_k^{m-1}u_m)^2}
$$
$$
+\sum\limits_{k=1}^n{(t_k-a)(v_1+\ldots+t_k^{m-1}v_m)^2}
+\sum\limits_{k=1}^n
{(b-t_k)(w_1+\ldots+t_k^{m-1}w_m)^2}
$$
$$
=\sum\limits_{i,j}{p_{i,j}u_iu_j} +\sum\limits_{i,j}{q_{i,j}v_iv_j}+\sum\limits_{i,j}{r_{i,j}w_iw_j}.
$$
It is easy to see that $p_{i,j}=s_{i+j-2}, \; q_{i,j}=s_{i+j-1}-as_{i+j-2}, \; r_{i,j}=bs_{i+j-2}-s_{i+j-1}$. Therefore,
 $H_m(a,b)$ is the matrix of the quadratic form $G$. If all $t_i\in [a,b]$, then  $G$ is nonnegative for any
${\bf u}, {\bf v}, {\bf w}$, i.e.  $G\succeq 0$.

Note that $r=\rk(R_n)$ is the number of distinct $t_i\; $ $=\; $ the total number of squares of linear forms (positive and negative) in the quadratic form $G({\bf u},{\bf 0},{\bf 0}).$ Then
$G({\bf u},{\bf 0},{\bf 0})\succeq 0$ yields that all $t_i$ are real numbers. Moreover,
$G({\bf 0},{\bf v},{\bf 0})\succeq 0$ and $G({\bf 0},{\bf 0},{\bf w})\succeq 0$ imply that all of them lie in $[a,b].$
\end{proof}

%If $m\ge n$, then  the Sylvester theorem shows that $t_1, \ldots, t_n$ belong to the interval $[a,b]$ if and only if $H_m(a,b)\succeq 0$. Moreover, the same holds for $m\ge k$=number of distinct $t_i$.

%Since all $t_i$ are real, we see that the number of distinct $\{t_i\}$ is equal to $\rk{H_{0,m}}.$

\medskip

In fact, the constraint $H_m(a,b)\succeq 0$ doesn't depend on $n=s_0$. Indeed, let
$$
\bar s_k:=\frac{t_1^k+\ldots+t_n^k}{n}=\frac{s_k}{n}, \quad \bar H_m(a,b)=\frac{1}{n} H_m(a,b).
$$
In other words, $ \bar H_m(a,b)$ can be obtained by substituting $\bar s_k$ for $s_k$ in
$H_m(a,b)$:
$$
\bar H_m(a,b)= \bar H_m(\bar s_1,\ldots, \bar s_{2m-1},[a,b]):=H_m(1,\bar s_1,\ldots,\bar s_{2m-1},[a,b]).
$$

Thus $H_m(a,b)\succeq 0$ if and only if $\bar H_m(a,b)\succeq 0$.

 Let us define in $(2m-1)$ - dimensional linear space $\{(\bar s_1,\ldots,\bar s_{2m-1})\}\cong {\bf R}^{2m-1}$ the domain
$$
\Delta_m(a,b):= \{(x_1,\ldots,x_{2m-1})\in {\bf R}^{2m-1}: \; \bar H_m(x_1,\ldots, x_{2m-1},[a,b])\succeq 0\}.
$$

\begin{theorem} For any real numbers $a, b, \; a<b, \; \Delta_m(a,b)$ is a compact convex  set in ${\bf R}^{2m-1}$ of dimensions $2m-1$.
\end{theorem}
\begin{proof} Since $-\infty<a<b<+\infty$ and $|x_k| \le
\max{(|a|^k,|b|^k)}$, we see that $\Delta_m(a,b)$ is a compact set. Note that $ H_m$ depends on $x_1, \ldots, x_{2m-1}$  linearly. Clearly, the convexity holds for any positive semidefinite matrix with entries that are linearly dependent on the variables  $x_1, \ldots, x_{2m-1}$. Since
$\bar H_m(a,b)\succ 0$ define a nonempty open set, we have $\dim(\Delta_m(a,b))=2m-1$.
%Let $C_+(a), \;
%C_-(b)\subset {\bf R}^{2m-1}$ are determined by
%$\bar H_m^+(a)\succeq 0, \; \bar H_m^-(b)\succeq 0$. It is not hard to prove that these domains are linearly equivalent to the convex cone
%$C_+(0)=-C_-(0)$. Thus, $\Delta_m(a,b)=C_+(a)\cap C_-(b)$.
\end{proof}

\section{$LP$ and $SDP_0$ bounds}

In this section we consider the linear programming  bound ($LP$) and the simplest version of the $SDP$ bounds  (we denote it by $SDP_0$)   for codes.

From here on we assume that $\Phi_k(t)$ is a polynomial of degree $k, \;
\Phi_k(\tau_0)=1$, and $S=T\cap [a,b]$ (the most interesting case for coding theory and sphere packings).

In fact, the optimal solution $A^*$ of the ${LPP}$ and  ${LPD}$ problems in Section 2 depends only on the family of polynomials $\Phi:=\{\Phi_k(t)\}_{k=0,1,\ldots}$ and $[a,b]$. We denote $1+A^*$ by
$LP(\Phi,[a,b])$.

\medskip

Since $\Phi_k(t)$ is a polynomial of degree $k$, we have:
$$\Phi_k(t)=p_{k0}+p_{k1}t+\ldots+p_{kk}t^k.$$
Then
$$
\Phi_k(t_1)+\ldots+\Phi_k(t_\ell)=\sum\limits_{d=0}^k{p_{kd}s_d}, \quad s_d=t_1^d+\ldots+t_\ell^d.
$$
%So for any $\Phi_k$ we have  uniquely defined coefficients $p_{k,0}, \ldots, p_{k,k}$.

Let $C=\{v_i\}$ be an $S$-code on ${\bf M}$, and let $\tau_{i,j}=\tau(v_i,v_j)$. Note that the number of ordered pairs $(v_i,v_j), \; i\ne j$, equals $n=|C|(|C|-1)$.  Then $(1)$ can be written in the form:
$$
y+p_{k0}+\sum\limits_{d=1}^k {p_{kd}x_d} \ge 0, \eqno (3)$$ where
$$
y=\frac{1}{|C|-1},  \quad
x_d=\bar s_d=\frac{s_d}{n}, \quad
s_d=\sum\limits_{i,j, i\ne j} {\tau_{i,j}^d}.
$$

From  Theorem 3.2  we have for any  $m$:
$$
\bar H_m(x_1,\ldots,x_{2m-1},[a,b])\succeq 0. \eqno (4)
$$

A positive  semidefinite constraints $(3)$ and $(4)$ make it possible to regard the problem of bounding $A({\bf M},S)$ as a semidefinite programming problem.
The standard form of the SDP problem is the following \cite{Todd,VB}:

\medskip

\noindent {\em Primal Problem}: $$ \mbox{minimize } \; \, c_1x_1+\ldots+c_{\ell}x_\ell$$
subject to
$$X\succeq 0, \; \mbox{ where }\; X=T_1x_1+\ldots+T_{\ell}x_{\ell}-T_0.$$

\noindent {\em Dual Problem}: $$ \mbox{maximize } \; \, \langle T_0,Y\rangle$$
subject to
$$\langle T_i,Y\rangle=c_i,\; i=1,\ldots {\ell},$$ $$ Y\succeq 0.$$

Here $T_0, T_1, \ldots, T_{\ell}, X,$ and $Y$ are real  $N\times N$ symmetric matrices,
$(c_1,\ldots,c_{\ell})$ is a cost vector, $(x_1, \ldots, x_{\ell})$ is a variable vector, and by
$\langle A,B\rangle$ we denote the inner product, i.e. $\langle A,B\rangle=\tr(AB)=\sum{a_{ij}b_{ij}}$.

Now we introduce the $SDP_0(\Phi,[a,b],m)$ bound.

\medskip

\noindent {$SDP_0$ {\em  Problem:}   Choose a natural number $m$ and real numbers $y, x_1, \ldots, x_{2m-1}$ so as to

$$\mbox{minimize } \; \, y$$
subject to
$$
y+\sum\limits_{i=1}^k {p_{ki}x_i} \ge -p_{k0}, \quad k=1,\ldots,2m-1,
$$
$$
\bar H_m(x_1,\ldots,x_{2m-1},[a,b])\succeq 0.
$$

\medskip

Note that in $(3)$: $|C|=(1+y)/y.$ Thus
\begin{theorem} If $y^*$ is the optimal solution of the $SDP_0$ problem, then
$$A({\bf M},S)\le SDP_0(\Phi,[a,b],m):=\frac {1+y^*}{y^*}.$$
\end{theorem}

%The $SDP_0$ problem looks very different from the $LP$ problem. However,
Since we just substituted
 $H_m\succeq 0$ for
$t \in S$ in the $LPP$ problem, we can expect that $\; SDP_0(\Phi,[a,b])=LP(\Phi,[a,b]).$  It is not  hard  to give a formal proof that $LPP=SDP_0$ and
$LPD=$ Dual $SDP_0$. We think that it becomes more clear by geometrical arguments. Actually, we are going to prove that $SDP_0=LPD$.

Consider in ${\bf R}^N,\; N=2m-1,$ the  domain defined by $(3)$:
$$
\Lambda_m(y):=\{{\bf x} \in {\bf R}^N: \sum\limits_{i=1}^k {p_{ki}x_i} \ge -y -p_{k0}, \; k=1,\ldots,N\}.
$$
Clearly, $\Lambda_m(y)$ is a  convex cone  in ${\bf R}^N$.

Let
$$
u_k=p_{k0}+\sum\limits_{i=1}^k {p_{ki}x_i}, \quad k=1,\ldots,N.
$$
%It follows from the orthogonality of $\{\Phi_k\}$
It is easy to see that $(u_1,\ldots,u_N)$
is a nondegenerate coordinate system of ${\bf R}^N$. So in these coordinates
$$
\Lambda_m(y)=\{(u_1,\ldots,u_N) \in {\bf R}^N: u_i \ge -y \; \mbox { for all } \; i \},
$$
i.e. $\Lambda_m(y)$ is a pointed cone with vertex $(-y,\ldots,-y)$.

For a given $m$ the
$SDP_0$ problem is asking to find the minimal $y=y^*(m)$ such that the intersection of $\Lambda_m(y)$ and $\Delta_m(a,b)$ is not empty. Then $y<y^*$ if and only if the intersection of these convex sets is empty.

For $0<y<y^*$ consider a separating hyperplane
$$
L({\bf u})=f_0+f_1u_1+\ldots+f_Nu_N
$$
such that
$$
L({\bf u})<0, \; \forall {\bf u}\in \Delta_m(a,b), \quad
L({\bf u})>0, \; \forall {\bf u}\in \Lambda_m(y).
$$
Clearly, $L({\bf u})>0$ yields $f_k\ge 0, f_0>y(f_1+\ldots+f_N)$. On the other hand, it can be proved that $L({\bf u})<0$ if and only if
$f(t)=f_0+f_1\Phi_1(t)+\ldots+f_N\Phi_N(t)<0$ for all $t\in [a,b]$. Then Theorem 2.4 yields $1+A^*=(1+y^*)/y^*$.

Thus we have proved:

\begin{theorem} $\; SDP_0(\Phi,[a,b])=LP(\Phi,[a,b]).$
\end{theorem}

In fact, for a continuous ${\bf M}$ the $LPP$ and $LPD$ problems are not finite linear programming problems.  These problems can be solved only via discretization. For instance, Odlyzko and Sloane \cite{OdS} ($\equiv$ \cite[Chapter 13]{CS}) applied $LPD$ for upper bounds on kissing numbers, where they replaced $S$ by 1001 equidistant points in $S$. For the $LPP$ problem it is not clear how to do a discretization of $\{\alpha_\tau\}$. 
On the other hand, for a given $m$ the $SDP_0$ is a finite primal SDP problem. As a by-product of  solutions of this problem we have bounds on $|C|$ and  power sums $s_k$
%That gives an interesting possibility of $SDP_0$'s  application  for constructions of optimal codes
(see the next section).

\medskip

\noindent{\bf Remark.} One of the important problems for $SDP_0$ is the problem of detecting an optimal $m$. We see that it is the same problem as the problem of optimal degree $N$ of a polynomial in Theorem 2.4. Our arguments show that $m=\lceil (N+1)/2 \rceil$.

The problem of finding optimal polynomials for Delsarte's linear programming bounds has been the subject of many papers. The most fundamental results in this problem were established by Levenshtein (see \cite{Levm}).

%Here we give several examples for  spherical ${\theta}$-codes on ${\bf S}^{n-1}$. If $\theta=\pi/2$, then $N=2$, i.e. $m=2.$ For $\theta=\pi/3$ (the kissing number problem) the degrees of suitable polynomials can be found in Table 1.5 \cite{CS}. For instance, if $n=8$, then $m=4$, and for $n=24$ we have $m=6$. However, $n=8, 24$ are exceptional dimensions. For other dimensions up to 23 we have: $m=5$ for $n=4$; $ m=6$ for $n=5,6,7,9,10,11,12$; $\; m=7$ for $n=13,\ldots, 21$; and $m=8$ for $n=22,23$.

%5. ``An extension of LP bound". Pfender's example +

\section{Bounds via $SDP_0$}

%Perhaps one of the most interesting  applications of the $SDP$ method is a possibility  for constructing  optimal (maximal) codes. Currently, we have just a draft of this approach,  the greater part of the work has yet to be done.
In this section we consider three examples of spherical codes and define the main steps of $SDP_0$.

\noindent{\bf 1.} It was shown by Levenshtein \cite{Lev2} and independently by Odlyzko and Sloane \cite{OdS}  that the maximal number of nonoverlapping unit spheres in ${\bf R}^8$ (resp. ${\bf R}^{24}$) that can touch another unit sphere is 240 (resp. 196560). In other words, they proved that the maximal size of a spherical $\pi/3$-code in dimension 8 (resp. 24) is 240 (resp. 196560). Bannai and Sloane \cite{BS}($\equiv$\cite[Chapter 14]{CS})  proved that the maximal arrangements in these dimensions are unique up to isometry.

Consider an application of the $SDP$ method for these spherical codes. Let us suppose  that we have no idea about the kissing number in eight dimensions and how a maximal kissing arrangement looks like. Now we try to find it by $SDP_0$.

\medskip

\noindent {\bf Step 1:}  {\em Semidefinite Programming}.

Let us apply the $SDP_0$ with $[a,b]=[-1,1/2]$ and $\Phi_k(t)=G_k^{(8)}(t)$. For $m=3, 4, 5$ we obtain that
$c^*(3)=324, \; c^*(4)= c^*(5)=240.$ (Here $c^*(m)=(1+y^*(m))/y^*(m)$ [see Section 4].) In fact, for any $m\ge 4$ we have $c^*(m)=240.$ So
we can expect that the maximal size of a kissing arrangement $C=\{v_i\}$ is 240 and for $m\ge 4$: $\rk(R_m)=4$  (the number of distinct $t_{ij}=\tau(v_i,v_j)$, see Section 3).  The $SDP_0$ for $m=4$  gives:
$$x_1^*=x_3^*=x_5^*=x_7^*=-1/239, \; x_2^*=29/239, \; x_4^*=8/239,\; x_6^*=11/956.$$
Thus if there exists $C$ with $c=|C|=240$ and $\rk(R_m)=4, \; m\ge 4$ , then
$$\sum \limits_{i=1}^4 \alpha_i\,t_i^k=\beta_k, \; k=0,\ldots,7, \; \mbox{ where } \; \beta_0=239, \; \beta_k:=239x_k^*, \; k\ge 1.$$

\medskip

\noindent {\bf Step 2 :} {\em Feasible solution for the distance distribution}.

Let us  consider the following problem:\\
Choose four real numbers $t_i, \; -1\le t_1\le t_2\le t_3\le t_4\le 1/2$ and four rational numbers
$\alpha_i,$ where $240\alpha_i$ are nonnegative integer numbers such that
%$$
%q_1+q_2+q_3+q_4=239;
%$$
$$
\alpha_1t_1^k+\alpha_2t_2^k+\alpha_3t_3^k+\alpha_4t_4^k=\beta_k,\quad k=0,\ldots, 7.
$$

Note that we have eight variables and eight equations. In fact,  this problem has the unique solution:

\begin{center}
\begin{tabular}{|c||c|c|c|c|}
%\hline
%$i$ & 1 & 2 & 3 & 4\\
\hline
$t_i$ & -1 & -1/2 & 0 & 1/2\\
\hline
$\alpha_i$ & 1 & 56 & 126 & 56\\
\hline
\end {tabular}
\end{center}
%$$
%t_1=-1, \; q_1=1, \; t_2=-1/2, \; q_2=56,\; t_3=0,\; q_3=126,\; t_4=1/2,\; q_4=56.
%$$
%(However, it is tricky to solve this problem algebraically, see  \cite{BS}.)

Geometrically we have the following picture: $\Delta_m(-1,1/2)\cap\Lambda_m(y)=\emptyset$ for $y<1/239$,
and for $y=1/239$ the intersection of these convex sets consists of  one point. In fact, this point is the vertex of the cone $\Lambda_m(y)$.

\medskip

\noindent {\bf Step 3 :} {\em Feasible solution for a maximal code}.

It is well known fact that $\{v_i\}$ are vectors in ${\bf R}^n$ if and only if its Gram matrix is positive semidefinite with rank that does not exceed $n$.

Let $C=\{v_1,\ldots,v_{240}\}$ be vectors in ${\bf R}^8$ and let $G=(t_{ij}),\; t_{ij}=\langle v_i,v_j\rangle$, be the Gram matrix of this set of vectors. We have ones on the diagonal of $G$ and other entries are $-1, 0, \pm 1/2.$  Let us denote by $Q_r(G)$ the total number of entries $r$ in $G$.
Now we have the following problem:\\
To find a symmetric $240\times 240$ positive semidefinite matrix $G$ with ones on the diagonal and with other entries $-1, 0, \pm 1/2$ such that $\rk(G)\le 8$ and
$$Q_{-1}(G)=240,\quad Q_{0}(G)=30240, \quad Q_{-1/2}(G)=Q_{1/2}(G)=13440.$$

Bannai and Sloane \cite [Theorem 7]{BS} show that this problem has the unique solution. Namely,
$C$ is isometric to the set of minimal vectors in $E_8$.

Note that if we apply this method for 24 dimensions, as well as for other spherical codes that were considered in \cite{BS}, then we also obtain optimal codes.

\medskip

\noindent{\bf 2.} Let us consider this approach for kissing arrangements in four dimensions. If we apply the $SDP_0$ for $m=3,4,5,6,7$, then we obtain: $$c^*(3)=c^*(4)=26,\quad c^*(5)=c^*(6)=c^*(7)\approx 25.5584.$$
From this it follows that the cardinality $c$ of $\pi/3$-spherical code $C$ is at most  $25.$

Let
$$\Omega_m(y):=\Delta_m(-1,1/2)\cap\Lambda_m(y).$$
As above, $\Omega_m(y^*(m)))$ consists of  one point. However, in four dimensions this point is not the vertex of $\Lambda_m(y)$.

%How can we find an optimal code in four dimensions?
%We know that $c\le 25.$ Suppose that $c=25$. Then for $m=5$ we have
%$\alpha_1+\ldots+\alpha_m=c-1=24$. On the other hand,  $25\alpha_i$  are nonnegative integer number. So we have just finite number possibilities for $\alpha_i$. Therefore, we can enumerate all possible cases for $\alpha_i$. If there exists a feasible solution for $\alpha_i, t_i$, i.e. $\{\alpha_i=\alpha_i^*\} \cap \Omega_m(y)\ne\emptyset$, then we go to Step 3. Otherwise, we suppose that $c=24$, and to repeat our computations for this case.

%In Step 3 for given $\{\alpha_i\}$ we try to find a feasible Gram matrix $G$. If there exists a feasible solution $G^*$, then we solved the problem. Otherwise, we consider other feasible $\{\alpha_i\}$ and if there are no feasible solutions for them, then we return to Step 2.
%(Actually, we know that there exists a solution for $c=24$ (see {\bf 3}).)

\medskip

\noindent{\bf 3.} If codes have some additional symmetries then  the problem of bounding of  codes sometimes can be solved easier than in general case.

Let us consider  spherical {\em antipodal} codes. We say that a spherical code $C$ is antipodal if $v\in C$ implies $-v\in C$. Note that the kissing number problem for antipodal codes in dimensions  $n\le 8$  was solved by Boyvalenkov \cite{Boyv}.

 It is easy to see that for an antipodal code $C$ we have
$$x_{2i+1}=-y=\frac{1}{1-c}.$$
So in this case the independent variables are $y$ and $x_2, x_4,\ldots,x_{2m-2}$.

If we apply the $SDP_0$ in four dimensions with $m=4$, then we obtain
$c^*=24, \; x_2^*=5/23, \; x_4^*=2/23, \; x_6^*=5/92$, and $x^*_{2i+1}=-1/23$. That yields (by the same method as in Step 2):
%
%\medskip
%
\begin{center}
\begin{tabular}{|c||c|c|c|c|}
\hline
%$i$ & 1 & 2 & 3 & 4\\
%\hline
$t_i$ & -1 & -1/2 & 0 & 1/2\\
\hline
$\alpha_i$ & 1 & 8 & 6 & 8\\
\hline
\end {tabular}
\end{center}
%
%\medskip
%$$
%t_1=-1, \; q_1=1, \; t_2=-1/2, \; q_2=8,\; t_3=0,\; q_3=6,\; t_4=1/2,\; q_4=8.
%$$
In fact, this implies the unique solution for the maximal code $C$: $C$ is isometric to the set of minimal vectors in $D_4$.

In the same way  maximal antipodal codes in dimensions 6,7 can be found. For $n=5$ the method gives
$c^*=42$. However, we expect that $c=40$. It is possible using the fact that $c\,\alpha_i$ are integer numbers to prove that $c\le 40$ (see details in \cite{Boyv}).

\section{An extension of Delsarte's method}

In this section we show how some recent extensions of Delsarte's method can be reformulated as a semidefinite programming problem.

Let $C=\{z_i\}$ be a finite subset of a compact 2-point-homogeneous space ${\bf M}$, and let
$$
s_k(C,[u,v]):=\sum\limits_{u\le\tau_{i,j}\le v}{\tau_{i,j}^k}=|C|\sum\limits_{u\le t\le v}{\alpha_t\,t^k},
\quad \tau_{i,j}=\tau(z_i,z_j);
$$
$$
s_k(C,[u,v)):=\sum\limits_{u\le\tau_{i,j}<v}{\tau_{i,j}^k}.
$$

The constraints $(3), (4)$ are relatively strong and sometimes give tight upper bounds for codes by Delsarte's linear programming method.  However, there are some additional constraints for $S$-codes which are not follow from $(3), (4)$.

For instance, Odlyzko and Sloane \cite{OdS} in dimension 17 made use of the following additional inequalities for a spherical $\pi/3$-code (a kissing arrangement):
$$
 \sum\limits_{t<t_1}{\alpha_t}\le 1, \quad   \sum\limits_{t<t_2}{\alpha_t}\le 2, \quad
t_1=-\frac{\sqrt{3}}{2}, \quad t_2=-\sqrt{2/3},
$$
i.e. $$s_0(C,[-1,t_1))\le |C|, \quad s_0(C,[-1,t_2))\le 2|C|.$$
These  inequalities easily can be derived from the following fact (see \cite{AVZ, BM} or \cite{Ran}):
the cardinality of a spherical $\theta$-code  on an open spherical cap  with angular radius $\varphi_k$ is at most $k$, where
$$\cos{\varphi_k}=t_k:=\sqrt{\cos{\theta}+(1-\cos{\theta})/(k+1)}.$$
This fact for a spherical $\theta$-code $C$  yields
$$
s_0(C,[-1,t_k))\le k\,|C|.
$$

Actually, any upper bound on the cardinality of a spherical $\theta$-code  in a spherical cap or  a    spherical strip implies some inequality for $s_0(C,[u,v])$. In \cite{BM} an overview   is given of the known bounds and  several new upper bounds for codes in spherical caps and strips have been proposed. Agrell, Vardy, and Zeger \cite{AVZ} applied  bounds on $s_0(C,[u,v])$ to improve some upper bounds for constant-weight codes.

Recently, Pfender \cite{Pf} found a new inequality for $s_2$:
% $C=\{x_i\}$ on a spherical cap of center $e_0$ and radius $\cos^{-1}{\sqrt{\cos{\theta}}}$ he proved that
$$
s_2(C,[-1,-\sqrt{\cos{\theta}}))\le s_0(C,[-1,-\sqrt{\cos{\theta}}))\cos{\theta}+ |C|(1-\cos{\theta}).
$$
Using this, he improved the upper bounds for the kissing numbers in dimensions 9, 10, 16, 17, 25 and 26. Note that Krasikov and Litsyn \cite{KrLi} using a similar approach improved some upper bounds for binary codes.

In the papers \cite{Mus, Mus13, Mus2, Mus3} we found a few inequalities for some linear combinations of $s_k(C,[-1,t_0])$ ($0\le k \le 9$), where $C$ is a spherical $\pi/3$-code and $t_0<-1/2.$
% on a spherical cap with radius $\varphi<\pi/3$.
In particular, that allowed us to prove that the kissing number in four dimensions is 24. This also yields a new solution of the Newton - Gregory (thirteen spheres) problem \cite{Mus13}.

Let $C$ be an $S$-code on  ${\bf M}$.  As above, we assume that
$S=T\cap[a,b]$ and $\Phi_k(t)$ is a polynomial of degree $k.$

Let $a=u_1<u_2<\ldots<u_\ell<u_{\ell+1}=b$ and let for $k=0,1,2\ldots$
$$
x_i^{(k)}:=s_k(C,[u_i,u_{i+1})), \quad i=1,\ldots,\ell-1; \quad x_\ell^{(k)}:=s_k(C,[u_\ell,b]).
$$
Then we have
$$
z=c+x_1^{(0)}+\ldots+x_\ell^{(0)}, \; \mbox{ where } \; z:=|C|^2,\; c:=|C|.
$$

In fact, the equality $z=c^2$ cannot be written as some positive semidefinite constraints. However,
the inequality $z \ge c^2$ can be expressed in the form:
$$\Gamma(c,z):=
\left(
\begin{array}{cccc}
 1 &  c\\
 c &  z
\end{array}
\right)
\succeq 0. \eqno (5)
$$

Using these notations $(3)$ can be expressed as:
$$
c+\sum\limits_{d=0}^k{p_{kd}}\sum\limits_{i=1}^\ell {x_i^{(d)}} \ge 0,\quad k=1,2,\ldots \eqno (6)
$$

Since $x_i^{(k)}$ are power sums of distances $\tau_{q,r}$ in the interval $[u_i,u_{i+1}]$,  we can apply Theorem 3.2 to this interval:
$$
H_m(x_i^{(0)},x_i^{(1)},\ldots,x_i^{(2m-1)},[u_i,u_{i+1}])\succeq 0, \quad i=1,\ldots,\ell. \eqno (7)
$$

Note that all additional inequalities that were considered above can be expressed in  $c, \{x_i^{(k)}\}$.   For instance, let $$\ell=2,\quad u_1=a=-1,\quad u_2=-\sqrt{\cos{\theta}},\quad u_3=b={\cos{\theta}}.$$ Then  Pfender's inequality can be expressed in the form:
$$
x_1^{(2)}\le x_1^{(0)}\,\cos{\theta}+ (1-\cos{\theta})\,c.
$$
This inequality is linear in $c, \{x_1^{(k)}\}$, i.e. it is a positive semidefinite constraint.

Suppose for a given $\ell$ and $\{u_i\}$ we have some additional positive semidefinite constraints for $c$ and $\{x_i^{(k)}\}$. Then we say that the set of these constraints is the set of $Ad$-constraints and write
$$
Ad(c,\{x_i^{(k)}\})\succeq 0. \eqno (8)
$$

\medskip

Using $(5)-(8)$ we introduce the following SDP problem:

\medskip

\noindent{\em $SDPA$ Problem:} For given $\ell, \{u_i\}$, $Ad$-constraints, and $m$ choose real numbers $c$ and $\; \left\{x_{i}^{(k)}\right\}_{i=1,\ldots,\ell}^{k=0,\ldots,2m-1}\; $ so as to
$$
\mbox{maximize } \; c
$$
subject to  $(5)-(8)$.
%$$ \Gamma(c,c+x_1^{(0)}+\ldots+x_\ell^{(0)})\succeq 0 \quad (\mbox{or equivalently} \; z=c+x_1^{(0)}+\ldots+x_\ell^{(0)}\ge c^2);$$
%$$
%c+\sum\limits_{d=0}^k{p_{k,d}}\sum\limits_{i=1}^\ell {x_i^{(d)}} \ge 0,\quad k=1,2,\ldots,2m-1;
%$$
%$$
%H_m(x_i^{(0)},x_i^{(1)},\ldots,x_i^{(2m-1)},[u_i,u_{i+1}])\succeq 0, \quad i=1,\ldots,\ell;
%$$
%$$
%Ad(c,\{x_i^{(k)}\})\succeq 0.
%$$

\medskip

%We denote the bound for this problem by
%$SDPA({\bf M},{\Phi},[a,b])$.

Denote by $c^*(m)$ the maximal value of $c$ that can be attained in this problem. If $C$ is an $S$-code then $c$ and
$x_i^{(k)}$ certainly satisfy the constraints in the $SDPA$ problem. Therefore $c^*(m)$ is an upper bound for $A({\bf M},S)$.
\begin{theorem} For any natural number $m$ we have
$$A({\bf M},S)\le c^*(m).$$
\end{theorem}

It is easy to see that $c^*(m+1)\le c^*(m).$ Let us denote the minimum of $c^*(m)$
 by $SDPA({\bf M},{\Phi},[a,b],m,Ad)$.

\begin{theorem}  $SDPA({\bf M},{\Phi},[a,b],m,Ad)\le SDP_0({\bf M},{\Phi},[a,b],m).$
\end{theorem}
\begin{proof}
Note that the variable $y$ in the $SDP_0$ problem can be expressed as
$$ y=\frac{c}{z-c}.$$
If we consider the constraints $(6), (7)$ in the $SDPA$ problem for the variables $y, \{x_i\}$ of the $SDP_0$ problem, then we obtain the constraints $(3), (4)$.
Therefore, the constraints in the $SDPA$ problem imply the constraints in the $SDP_0$ problem.
(Moreover, we have the additional constraint $(8).$)
Then we have
$$ \frac{c^*}{z-c^*}\ge y^*.$$
We also have $(5)$: $z\ge c^2$ which yields
$$
 \frac{1}{c^*-1}=\frac{c^*}{(c^*)^2-c^*}\ge\frac{c^*}{z-c^*}\ge y^*.
$$
Thus
$$ c^*\le \frac{1+y^*}{y^*}. $$
\end{proof}

Clearly that the $SDPA$ problem with empty $Ad$-constraints set has the same constraints as the $SDP_0$ problem. Then $$SDPA({\bf M},{\Phi},[a,b],m,\emptyset) = SDP_0({\bf M},{\Phi},[a,b],m),$$ i.e. the $SDPA$ bound is equal to the Delsarte linear programming bound. Therefore, the $SDPA$ bound  with nonempty $Ad$-constraints set can be considered as {\em an extension of Delsarte's method}. We see that for some cases this extension gives better bounds than the classical Delsarte method
\cite{KrLi, Mus, Mus13, Mus2, Mus3, Pf}. However, it is not easy to find (and to prove) new $Ad$-constraints.

\section{ $SDP$ bounds for codes }

%The constraints $(3)$  in the $SDP_0$ problem are linear.
%These constraints derive from Corollary 2.1.
% For a linear programming problem we have to have only linear inequalities, but for an SDP problem it is not necessarily. This method allows positive semidefinite constraints.  Note that $(3)$ is a simple consequence of the p.d.k (positive semidefinite) property given in Theorem 2.1.
In this section we explain how to use the p.d.k. property to obtain upper bounds on $S$-codes in a subset
${\bf \Pi}\subseteq {\bf M}$ and to improve the $SDP_0$ and $SDPA$ bounds.

For any finite subsets  $X, Y$ in a 2-point-homogeneous space ${\bf M}$ we define by
$$
F_k(X,Y):=\sum\limits_{x\in X, y\in Y} {\Phi_k(\tau(x,y))}.
$$
It is easy to see that Theorem 2.1 implies the following fact:

\begin{theorem}
Let $\{X_i\}_{i=1,\ldots d}$ be finite subsets of ${\bf M}.$ Then for any
%finite $X\subset{\bf M}$ and
 integer $k\ge 0$ the matrix $\left(F_k(X_i,X_j)\right)$ is positive semidefinite.
\end{theorem}

%Note that the case ${\bf\Pi}\subset{\bf M}, \; {\bf\Pi}\ne{\bf M}$ is also interesting for applications.
%For instance,
%constant-weight codes \cite{CS} and codes in spherical caps \cite{AVZ, BM} are extensively studied in coding theory.
%It is easy to extend the $SDP_0$ bounds (as well as the LP bounds)  for codes on ${\bf\Pi}$. However, we are not sure that it gives sharp bounds for  ${\bf\Pi}\ne {\bf M}$.
%In this section we consider an extension of the $SDP_0$ bounds for codes on ${\bf\Pi}\subseteq{\bf M}$.
%However, if ${\bf\Pi}\ne{\bf M}$, then $SDPA$ can give much better upper bounds on $A({\bf\Pi},S)$ than $SDP_0$.
%For instance, let ${\bf\Pi}$ be a spherical cap in ${\bf M}={\bf S}^{n-1}$ of center $e_0$ and radius $\varphi.$ In this case:
%$a_{01}=\cos{\varphi}, \; b_{01}=1.$ Then for a ${\theta}$-code on ${\bf\Pi}$ with $\theta>\varphi$ the $SDPA$ bound  certainly  is less than the $SDP_0$ bound. Moreover, for some cases the $SDPA$ bounds are sharp. The same holds for constant weight codes.
%(Of course,
%${\bf\Pi}={\bf M}$  is the most interesting case for us.)

%Let $e_0$ be a fixed point in ${\bf M}$.  In some cases $e_0$ is uniquely defined by the problem. For instance, if ${\bf\Pi}$ is a spherical cap, then $e_0$ is the center of this cap. In the case
%${\bf\Pi}={\bf M}$ a good choice for $e_0$ is when $e_0$ is a point of an $S$-code.

%\medskip

%Now we replace the linear inequalities $(3)$ by  semidefinite constraints.

%Consider
%$S$-codes on  ${\bf\Pi}$.
%Let $T({\bf\Pi}):=\{\tau(x,y): x, y\in {\bf\Pi}\}.$

\medskip

\noindent{\bf 7-A.  The $SDPA$ bounds for codes on subsets.}
Let ${\bf\Pi}$ be any subset of ${\bf M}$. The largest size of an $S$-code in ${\bf\Pi}$ will be denoted by $A({\bf\Pi},S)$. (As above, we assume that
$S=T\cap[a,b].$)

Let $C_0=\{q_1, \ldots, q_r\}$ be a fixed subset of  ${\bf M}$. Let $-\infty\le\alpha_i\le\beta_i\le\infty$ for $i=1,\ldots,r$. We define ${\Pi_1}={\Pi_1}(\{\alpha_i,\beta_i\})$  by the following constraints:
$$
{\Pi_1}=\{x\in{\bf M}: \alpha_i\le\tau(x,q_i)\le\beta_i\}, \quad i=1,\ldots,r.
$$
 Clearly,  if $\alpha_i=-\infty,\; \beta_i=\infty$ for all $i$, then  ${\Pi_1}={\bf M}.$

  Suppose ${\bf\Pi}\subseteq{\Pi_1}$
Consider any $S$-code $C$ on ${\bf\Pi}$. Let $X_i=\{q_i\}, \; i=1,\ldots,r,$ and $X_{r+1}=C$. Then Theorem 6.1 yields
$$ G_k:=(F_k(X_i,X_j))\succeq 0, \quad k=0, 1,\ldots . \eqno (9)$$

Note that in the matrix $G_k$ for $1\le i,j\le r$ entries $(G_k)_{ij}=\Phi_k(\tau(q_i,q_j))$  are given.   $(G_k)_{r+1,r+1}=F_k(C,C)$ (see $(6)$).
Let
$$y_i^{(k)}:=\sum\limits_{q\in C}{\tau^k(q,q_i)}.$$
Then $y_1^{(0)}=y_2^{(0)}=\ldots= y_r^{(0)}=c=|C|$ and
$$(G_k)_{i,r+1}=\sum\limits_{d=0}^k {p_{kd}\,y_i^{(d)}}, \quad i\le r.$$
 Thus, the matrix $G_k$ linearly depends on the variables $y_i^{(k)}$ and $x_i^{(k)}$.

On the other hand, since $\alpha_i\le\tau(q,q_i)\le\beta_i, \; q\in C,$ we have
$$
H_m(y_i^{(0)},y_i^{(1)},\ldots,y_i^{(2m-1)},[\alpha_i,\beta_{i}])\succeq 0, \quad i=1,\ldots,r. \eqno (10)
$$

For a given $m$ let us consider the following SDP problem:
$$\mbox{maximize } \; c$$
s.t.
$$ (7),\; (8),\; (9),\; (10).$$

Denote by $c^*(m)=c^*(m,\{q_i\},\{\alpha_i\}, \{\beta_i\})$ the maximal value of $c$ that can be attained in this problem, and let us denote the minimum of $c^*(m)$ by $SDPA({ \Pi_1},{\Phi},[a,b],Ad)$. Then we have the following theorem:

\begin{theorem} $A({\bf\Pi},S)\le SDPA({ \Pi_1},{\Phi},[a,b],Ad)\le c^*(m)$
\end{theorem}

%
%Let ${\bf \Pi}\subseteq {\bf M}$ and $e_0\in {\bf M}$.

\medskip

\noindent{\bf 7-B.  The $\widehat{SDP}$ bounds.} Now using Theorem 7.1 we improve the $SDPA$ bound.
Let $C_0=\{q_1,\ldots,q_r\}$ be a fixed subset of ${\bf M}$ and let ${\bf\Pi}\subseteq{\bf M}$. Suppose  a family $\{\Pi_i\}_{i=1,\dots,d}$  is a partitioning  of ${\bf\Pi}$,
$${\bf\Pi}\subseteq\Pi_1\cup\ldots\cup \Pi_d, \quad \Pi_i\cap\Pi_j=\emptyset, \; i\ne j.$$
We also assume that
$$\Pi_j\subseteq \bar\Pi_j:=\{q\in{\bf M}: \alpha_{ij}\le\tau(q_i,q)\le\beta_{ij}, \; i=1,\ldots, r\}.$$
%However, for applications be better to use  $\Pi_i$ with this property.
For finite  $C\subset{\bf\Pi}$ denote  by
$$C_i:=\Pi_i\cap  C,  \quad
 c_i:=|C_i|. $$
%We have $c_0=1.$
%Now we apply theorem 5.1 for $\; \tilde C=\{e_0\}\cup C$,
%Let $C$ be a finite subset in ${\bf\Pi}$, and let
$$
x^{(k)}_{ij}:=\sum\limits_{u\in C_i}\sum\limits_{v\in C_j}{\tau^k(u,v)}, \quad 1\le i, j \le d, \quad
k=0,1,2,\ldots,
$$
$$
y^{(k)}_{ij}:=\sum\limits_{q\in C_j}{\tau^k(q_i,q)}, \quad 1\le i\le r, \quad 1\le j \le d, \quad
k=0,1,2,\ldots .
$$

Recall that $\Phi_k(t)=\sum {p_{ki}t^i}.$  Then
$$
F_k(C_i,C_j)=\sum\limits_{\ell=0}^k {p_{k\ell}\,x^{(\ell)}_{ij}},
\quad F_k(q_i,C_j)=\sum\limits_{\ell=0}^k {p_{k\ell}\,y^{(\ell)}_{ij}}.
$$

If we apply Theorem 7.1 with $X_i=\{q_i\},  i=1,\ldots, r, \; X_{i+r}=C_i, i=1,\ldots,d$, then we obtain that the matrix
$\Gamma_k=(F_k(X_i,X_j)$ is positive semidefinite. We see that entries of $\Gamma_k$ are linearly dependent on the variables ${\bf x}=\left\{x^{(k)}_{ij}\right\}$ and ${\bf y}=\left\{y^{(k)}_{ij}\right\}$. (Obviously,
$x^{(k)}_{ji}=x^{(k)}_{ij}$. So we can assume that ${\bf x}$ is a set of variables $x^{(k)}_{ij}$ with $i\le j$.)
Therefore, $\Gamma_k=\Gamma_k({\bf x},{\bf y}).$
 Thus for any $m>0$: %we have
$$ \Gamma_k({\bf x},{\bf y}) \succeq 0, \quad k=0,1,\ldots,2m-1. \eqno (11)$$

%Now we take a look at ${\bf x}^{(0)}$ as at the set of variables. We obviously have
Note that all $x^{(0)}_{ij}, y^{(0)}_{ij}$ are {\em nonnegative integer} numbers and
$$
 y^{(0)}_{ij}=c_j,  \quad x^{(0)}_{ij}=c_ic_j. \eqno (12)
$$

Since $\alpha_{ij}\le\tau(q,q_i)\le\beta_{ij}, \; q\in C_j,$ we have
$$
H_m(y_{ij}^{(0)},y_{ij}^{(1)},\ldots,y_{ij}^{(2m-1)},[\alpha_{ij},\beta_{ij}])\succeq 0; \quad i=1,\ldots,r;\; j=1,\ldots,d. \eqno (13)
$$

\medskip

Let $C$ be an $S$-code  on ${\bf\Pi}$. All  values $\{\tau(u,v)\}, u\in C_i, v\in C_j, u\ne v$ belong to some interval $[a_{ij},b_{ij}]\subseteq [a,b]$. It is clear that the bounds
$\{a_{ij}, b_{ij}\}$ depend only on $\Pi_i, \Pi_j$, and $[a,b]$.
%Indeed, for $0<i\le j\le d$ we have
%$$
%[a_{ij},b_{ij}]=[\alpha_{ij},\beta_{ij}]\cap [a,b],$$
%where
%$$
%\alpha_{ij}=\inf\limits_{u\in\Pi_i,v\in\Pi_j} {\tau(u,v)}, \quad
%\beta_{ij}=\sup\limits_{u\in\Pi_i,v\in\Pi_j} {\tau(u,v)},
%$$
%and for $i=0$ we have
%$$
%a_{0j}=\inf\limits_{u\in\Pi_j} {\tau(e_0,u)}, \quad b_{0j}=\sup\limits_{u\in\Pi_j} {\tau(e_0,u)}.
%$$
%Indeed,
%yields $c_i^2\le y_{ii},\; y_{ij}^2\le y_{ii}\,y_{jj}$, etc. It is not hard to prove for $d=1$ that if we consider $SDP_0$ in these variables, then  maximizing $c_1$ we obtain the same bound.
%Now we apply Theorem 3.2 for three cases: (i) $i\ne j$, (ii) $i=j$, (iii) $i=0.\\$
%(i).

By assumption for $i\ne j$ we have $\Pi_i\cap\Pi_j=\emptyset$ which yields $C_i\cap C_j=\emptyset$. Then $\tau(u,v)\in [a_{ij},b_{ij}]$ for all  $u\in C_i, v\in C_j$. Therefore,
%we have:
$$
H_m(x^{(0)}_{ij},x^{(1)}_{ij},\ldots,x^{(2m-1)}_{ij},[a_{ij},b_{ij}])\succeq 0, \quad 1\le i<j \le d.
\eqno (14)
$$

Note that among ordered pairs $(u,v), \; u, v\in C_i$ we have exactly  $c_i$ pairs $(u,u)$ such that $\tau(u,u)=\tau_0\notin [a_{ii},b_{ii}].$ All other pairs belong to this interval. Let $z_{k,i}:=x_{ii}^{(k)}-c_i\,\tau_0^k,\; k>0,\; \,
z_{0,i}:= x_{ii}^{(0)}-c_i$. Then we have
$$
H_m(z_{0,i},z_{1,i},\ldots,z_{2m-1,i},[a_{ii},b_{ii}])\succeq 0,\quad i=1,\ldots, d. \eqno (15)
$$
%(iii) In this case
%and
%$$
%H_{m}^{0,j}({\bf x},[a_{0j},b_{0j}]):=H_m(x^{(0)}_{0j},x^{(1)}_{0j},\ldots,x^{(2m-1)}_{0j},[a_{0j},b_{0j}])\succeq 0, \quad j=0,\ldots,d.
%$$
%Note that the case $(i,j)=(0,0)$ is trivial. Finally,
%$$H_{m}^{i,j}({\bf x},[a_{ij},b_{ij}])\succeq 0 \eqno (10)$$
%is well defined for all $i, j, m$, where $0\le i,j\le d,\; m=0,1,2\ldots .$
%\medskip

In Section 6 we considered some additional constraints which don't follow from $(3)$, $(4)$.
%In fact, these constraints are linear inequalities in the variables ${\bf x}$, i.e. they are the SDP constraints.
In the case when  with $(11), (13), (14), (15)$ we have some additional positive semidefinite constraints
%for ${\bf x}$ and ${\bf y}$
we write them in the following form:
%we say that the set of these constraints is an  $Ad$ - constraints set  and write
$$ Ad({\bf x},{\bf y},\{\Pi_i\},[a,b])\succeq 0. \eqno (16) $$

\medskip

Now we introduce the $QOP$ bounds for codes on ${\bf\Pi}$.

\medskip

%Input:\\
%Constants: $\; {\bf p}=\{p_{k,i}\}, a, b$\\
%Fixed parameters: $%\Phi=\{\Phi_k(t)\}_{k=0,1,\ldots}, \; \,
%d, \, \{\Pi_i\}$ ($a_{ij}, b_{ij}$ are uniquely defined by these parameters)\\
%(i.e. $\; \left\{[a_{ij},b_{ij}]\right\}_{\forall i,j:0\le i\le j\le d}$) \\
% {\bf\Pi}=\{\Pi_i\}_{i=0,1,\ldots,d}
%Variable:  $m$

%Output: $c^*(m)$

%\medskip

\noindent{\em $QOP$ Problem:} For given $m, C_0, d, \{\Pi_i\}$, and $Ad$-constraints  choose nonnegative integer numbers $c_i$ and real numbers ${\bf x}, {\bf y}$
%real numbers $\; \left\{x_{ij}^{(k)}\right\}_{\forall i,j:0\le i\le j\le d}^{k=1,\ldots,2m-1}\; $ with  $x_{00}^{(k)}=\tau_0^k,$ and
% $\; \left\{x_{ij}^{(0)}\right\}\; $ such that
%$\; x_{00}^{(0)}=1, \; \, x^{(0)}_{ij}=x^{(0)}_{0i}\,x^{(0)}_{0j} \; $
so as to
$$
\mbox{maximize } \; c=c_1+\ldots+c_d
$$
subject to
$$ (11)-(16).$$

\medskip

%$$ \Gamma_k({\bf x},{\bf p}) \succeq 0, \quad k=1,\ldots,2m-1;$$
%$$H_{m}^{i,j}({\bf x},[a_{ij},b_{ij}])\succeq 0, \quad \forall i,j: 0\le i\le j\le d;$$
%$$ Ad({\bf x},\{\Pi_i\},[a,b])\succeq 0.$$
%Since $H_m\succeq 0$ implies $H_{m-1}\succeq 0$, we have $c^*(m)\le c^*(m-1).$ Denote by $c^*$ the minimal value $c^*(m)$.
%The $QOP$ bounds for codes depend on $d$ and a family $\{\Pi_i\}$. Clearly,  if we denote for given $d$ and an optimal $\{\Pi_i\}$ the bound $c^*$ by $QOP_d({\bf\Pi},{\Phi},[a,b])$, then
%  $$QOP({\bf\Pi},{\Phi},[a,b])\le QOP_{d+1}({\bf\Pi},{\Phi},[a,b])\le QOP_d({\bf\Pi},{\Phi},[a,b]).$$
%Finally, we have
%$$
% A({\bf \Pi},S)\le QOP({\bf\Pi},{\Phi},[a,b])\le SDPA({\bf\Pi},{\Phi},[a,b])\le LP({\bf\Pi},{\Phi},[a,b]).
%$$
%It is extremely important open question for this problem: how to define suitable $d, \; \{\Pi_i\}$, and $m$? In the next section we will consider several examples for spherical codes.

Actually, the $QOP$ problem is not a semidefinite programming problem.
In $(12)$ we have $(i)$ $x_{ij}^{(0)}=c_ic_j$, and $(ii)$ $c_i$ are nonnegative {\em integer} numbers. So $(i)$ and $(ii)$ are not positive semidefinite constraints.
%Another interesting question is how  using this fact to improve the $SDP$ bounds?

%Currently we have no  efficient algorithms for  the $QOP$ problem. One of the natural ideas is to replace $(i)$ by some positive semidefinite constraints and to solve the problem using SDP perhaps with some integer approximations based on $(ii)$.

Let us consider  $c_i$ and $x_{ij}^{(0)}$ as independent real variables
%But in this case we are loosing the linearity in $\Gamma_k$, and therefore we cannot apply $(5)$ for SDP. However,
with the following positive semidefinite constraints:
$$  y_{ij}^{(0)}=c_j\ge 0, \quad
%\Gamma_0=
\left(
\begin{array}{cccc}
 1 &  c_1
&\ldots &  c_d\\
 c_1 &  h_{11}
&\ldots &  h_{1d}\\
\vdots & \vdots
&\ddots & \vdots\\
 c_d &  h_{1d}
&\ldots &  h_{dd}
\end{array}
\right)
\succeq 0,
\quad h_{ij}=x_{ij}^{(0)}. \eqno (17)
$$
If we substitute $(17)$  with real $c_i$  for $(12)$ in the $QOP$ problem, then we get a semidefinite programming problem. Let us denote this problem by $\widehat{SDP}$.
%({\bf\Pi},{\Phi},[a,b])$

Denote by $c^*(m,C_0,d,\{\Pi_i\},Ad)$ the maximal value of the sum $c$ that can be attained in the $\widehat{SDP}$ problem. Then
%If $C$ is an $S$-code then
%${\bf x}$ certainly satisfies the constraints in the $QOP$ problem. Therefore $c^*(m)$ is an upper bound for $A({\bf \Pi},S)$.
%solution of this problem is an upper bound for $A({\bf M},S)$.
%Thus we have proved:
%Suppose a family  of subsets $\{\Pi_i\}$ is such that  ${\bf\Pi}=\Pi_1\cup\ldots\cup \Pi_d$, and
%$\Pi_i\cap\Pi_j=\emptyset$ for $i\ne j$. Then for any natural number $m$ we have
\begin{theorem}
$A({\bf\Pi},S)\le c^*(m,C_0,d,\{\Pi_i\},Ad).$
\end{theorem}

\medskip

\noindent{\bf 7-C. New upper bounds for one-sided kissing numbers.}
Our first computer experiments show that this method can significantly to improve some upper bounds for codes.
% The main question is: how good the $\widehat{SDP}$ bounds are? Could these bounds essentially to improve the $LP$ bounds?  that with  the $\widehat{SDP}$ bounds just slightly better than $LP$. In these experiments we are working only with spherical codes, where
%$\Pi_i$ are spherical strips. Perhaps
To show that we consider the one-sided kissing problem.

Let $H$ be a closed half-space of ${\bf R}^{n}$. Suppose $S$ is a unit sphere in $H$ that touches the supporting hyperplane of $H$. The one-sided kissing number $B(n)$ is the maximal number of unit nonoverlapping spheres in $H$ that can touch $S$.

If  nonoverlapping  unit spheres kiss (touch) the unit sphere $S$ in $H\subset{\bf R}^n$, then the set of kissing points
is an arrangement on the closed hemisphere $S_+$ of $S$ such that the (Euclidean) distance between any two points is at least 1. So the one-sided kissing number problem can be stated in another way: How many points can be placed on the surface of $S_+$  so that the angular separation between any two points is at least $\pi/3$? In other words, $B(n)$ is the maximal cardinality of a $\pi/3$-code on the hemisphere $S_+$.

Currently $B(n)$ known only for $n\le 4$. Clearly, $B(2)=4$. It is not hard to prove that $B(3)=9.$
 Recently, using some extensions of Delsarte's method we proved that $B(4)=18$ (see \cite{Mus3} for a proof and references).

In \cite[Section 7]{BM} we derive estimates of $B(n)$ for $n=5,6,7,8.$ Actually, these bounds are based on the LP bounds for kissing numbers and the LP bounds on $\theta_0$-codes, where
$\cos{\theta_0}=1/\sqrt{3}.$ So we call these upper bounds  the LP bounds.

We applied the $\widehat{SDP}$ with $[a,b]=[-1,1/2]$, the empty $Ad$ - constraints set, $m=6,\; r=1$  (i.e. $C_0$ consists of the one point $e_0$), $d=3,$ and
$$
\Pi_1=\{p\in {\bf S}^{n-1}: \dist(p,e_0)\le \pi/4\},
$$
$$
\Pi_2=\{p\in {\bf S}^{n-1}: \pi/4<\dist(p,e_0)< 2\pi/5\},
$$
$$
\Pi_3=\{p\in {\bf S}^{n-1}: 2\pi/5\le\dist(p,e_0)\le \pi/2\}.
$$
We give the $\widehat{SDP}$ (as well as LP) bounds on $B(n)$ in the following table:

\begin{center}
\begin{tabular}{|c||c|c|c|c|c|c|c|}
\hline
$n$ & 3 & 4 & 5 & 6 & 7 & 8 & 9 \\
\hline
Lower bound known & 9 & 18 & 32 & 51 & 93 & 183 & \\
\hline
LP bound & 9 & 20 & 39 & 75 & 135 & 238 & 378\\
\hline
SDP bound & 9 & 19 & $35$ & 64 & $110$ & $186$ & $309$\\
\hline
\end {tabular}
\end{center}
%
%\medskip
%$$
%t_1=-1, \; q_1=1, \; t_2=-1/2, \; q_2=8,\; t_3=0,\; q_3=6,\; t_4=1/2,\; q_4=8.
%$$

Note that for $n=8$ the $\widehat{SDP}$  bound is close to the known lower bound: $B(8)\ge 183.$
So using some new $Ad$ - constraints or other families $\{\Pi_i\}$  we have a good chance to solve the one-sided kissing problem in eight dimensions. \footnote{Recently, the equality $B(8)=183$ has been proved by Bachoc and Vallentin \cite{bac06b} via SDP. They are also  improved our bounds on $B(n)$ for $5\le n\le 10$.}

%\medskip

%\begin{theorem} Suppose\\
%(i) \, If $\; m\ge n$, then $c^*(m)\le c^*(n)$.\\
%(ii) If $\; {\widehat{\bf\Pi}}\sqsupseteq {\bf\Pi}$,  then $c^*(m,{\widehat{\bf\Pi}})\le c^*(m,{\bf\Pi})$.
%\end{theorem}
%\begin{proof}
%It is clear that $c^*(m)$ is not increasing whenever $m$ is increasing.
%\end{proof}

%\section{On SDP Bounds for Spherical Codes}

\end{document}